\documentclass{article}

\usepackage{latexsym,epsfig}
\usepackage{amssymb,amsmath,amsfonts}
\usepackage{exscale} 
\usepackage{ifthen}
\usepackage{longtable}
\usepackage{subfigure}
\usepackage{caption}
\usepackage{cite}
\usepackage{lscape}
\usepackage{alltt} 
\usepackage{multicol} 

\begin{document}
\bibliographystyle{gabialpha}
\protect\pagenumbering{arabic}
\setcounter{page}{1}
 
\newcommand{\Zt}{\rm}

\newcommand{\ba}{\begin{array}}
\newcommand{\ea}{\end{array}}
\newcommand{\pot}{{\cal P}}
\newcommand{\curv}{\cal C}
\newcommand{\ddt} {\mbox{$\frac{\partial  }{\partial t}$}}
\newcommand{\hl}{\sf}
\newcommand{\hd}{\sf}

\newcommand{\Ad}{\mbox{\rm Ad}}
\newcommand{\Adsm}{\mbox{{\rm \scriptsize Ad}}}
\newcommand{\ad}{\mbox{\rm ad}}
\newcommand{\adsm}{\mbox{{\rm \scriptsize ad}}}
\newcommand{\diag}{\mbox{\rm Diag}}
\newcommand{\sect}{\mbox{\rm sec}}
\newcommand{\id}{\mbox{\rm id}}
\newcommand{\idsm}{\mbox{{\rm \scriptsize id}}}
\newcommand{\eps}{\varepsilon}

\newcommand{\aL}{\mathfrak{a}}
\newcommand{\bL}{\mathfrak{b}}
\newcommand{\mL}{\mathfrak{m}}
\newcommand{\kL}{\mathfrak{k}}
\newcommand{\gL}{\mathfrak{g}}
\newcommand{\nL}{\mathfrak{n}}
\newcommand{\hL}{\mathfrak{h}}
\newcommand{\pL}{\mathfrak{p}}
\newcommand{\uL}{\mathfrak{u}}
\newcommand{\lL}{\mathfrak{l}}

\newcommand{\kG}{{\tt k}}
\newcommand{\nG}{{\tt n}}

\newcommand{\Cart}{$G=K e^{\overline{\aL^+}} K$}
\newcommand{\Area}{\mbox{Area}}
\newcommand{\Hd}{\mbox{\rm Hd}}
\newcommand{\Hdim}{\mbox{\rm dim}_{\mbox{\rm \scriptsize Hd}}}
\newcommand{\Tr}{\mbox{\rm Tr}}
\newcommand{\bs}{{\cal B}}
\newcommand{\nc}{{\cal N}}
\newcommand{\MM}{{\cal M}}
\newcommand{\Ch}{{\cal C}}
\newcommand{\clCh}{\overline{\cal C}}
\newcommand{\Cnt}{\mbox{\rm C}}

\newcommand{\NN}{\mathbb{N}} \newcommand{\ZZ}{\mathbb{Z}}
\newcommand{\QQ}{\mathbb{Q}} \newcommand{\RR}{\mathbb{R}}
\newcommand{\KK}{\mathbb{K}} \newcommand{\FF}{\mathbb{F}}
\newcommand{\CC}{\mathbb{C}} \newcommand{\EE}{\mathbb{E}}
\newcommand{\XX}{X}
\newcommand{\HH}{I\hspace{-2mm}H}
\newcommand{\norm}{\Vert\hspace{-0.35mm}|}
\newcommand{\Sph}{\mathbb{S}}
\newcommand{\ganz}{\overline{\XX}}
\newcommand{\rand}{\partial\XX}
\newcommand{\prodrand}{\partial\XX_1\times\partial\XX_2} 
\newcommand{\regrand}{\partial\XX^{reg}}
\newcommand{\singrand}{\partial\XX^{sing}}
\newcommand{\Frand}{\partial\XX^F}
\newcommand{\Lim}{L_\Gamma}          
\newcommand{\cLim}{M_\Gamma}          
\newcommand{\Flim}{K_\Gamma}
\newcommand{\reglim}{L_\Gamma^{reg}}
\newcommand{\radlim}{L_\Gamma^{rad}}
\newcommand{\raylim}{L_\Gamma^{ray}}
\newcommand{\horinf}{\mbox{Vis}^{\infty}}
\newcommand{\horF}{\mbox{Vis}^B}
\newcommand{\Sml}{\mbox{Small}}
\newcommand{\SmlF}{\mbox{Small}^B}

\newcommand{\ifl}{\qquad\Longleftrightarrow\qquad}
\newcommand{\at}{\!\cdot\!}
\newcommand{\ging}{\gamma\in\Gamma}
\newcommand{\xo}{{o}}
\newcommand{\gamo}{{\gamma\xo}}
\newcommand{\gam}{\gamma}
\newcommand{\gax}{h}
\newcommand{\gxi}{{G\!\cdot\!\xi}}
\newcommand{\bd}{$(b,\Gamma\at\xi)$-densit}
\newcommand{\bt}{$(b,\theta)$-densit}
\newcommand{\cd}{$(\alpha,\Gamma\at\xi)$-density}
\newcommand{\be}{\begin{eqnarray*}}
\newcommand{\ee}{\end{eqnarray*}}

\newcommand{\an}{\ \mbox{and}\ }
\newcommand{\as}{\ \mbox{as}\ }
\newcommand{\diam}{\mbox{diam}}
\newcommand{\is}{\mbox{Isom}}
\newcommand{\Ax}{\mbox{Ax}}
\newcommand{\Fix}{\mbox{Fix}}
\newcommand{\Par}{F}
\newcommand{\Min}{\mbox{Fix}}
\newcommand{\vol}{\mbox{vol}}
\newcommand{\Td}{\mbox{Td}}
\newcommand{\piF}{\pi^B}
\newcommand{\piKM}{\pi^I}

\newcommand{\for}{\ \mbox{for}\ }
\newcommand{\pr}{\mbox{pr}}
\newcommand{\sh}{\mbox{sh}}
\newcommand{\shi}{\mbox{sh}^{\infty}}
\newcommand{\rank}{\mbox{rank}}
\newcommand{\supp}{\mbox{supp}}
\newcommand{\mass}{\mbox{mass}}
\newcommand{\kernel}{\mbox{kernel}}
\newcommand{\st}{\mbox{such}\ \mbox{that}\ }
\newcommand{\Stab}{\mbox{Stab}}
\newcommand{\Root}{\Sigma}
\newcommand{\Cone}{\mbox{C}}
\newcommand{\wrt}{\mbox{with}\ \mbox{respect}\ \mbox{to}\ }
\newcommand{\where}{\ \mbox{where}\ }

\newcommand{\con}{{\sc Consequence}\newline}
\newcommand{\rem}{{\sc Remark}\newline}
\newcommand{\prf}{{\sl Proof.\  }}
\newcommand{\qed}{$\hfill\Box$}

\newenvironment{rmk} {\newline{\sc Remark.\ }}{}  
\newenvironment{rmke} {{\sc Remark.\ }}{}  
\newenvironment{rmks} {{\sc Remarks.\ }}{}  
\newenvironment{nt} {{\sc Notation}}{}  

\newtheorem{satz}{\bf Theorem}

\newtheorem{df}{\sc Definition}[section]
\newtheorem{cor}[df]{\sc Corollary}
\newtheorem{thr}[df]{\bf Theorem}
\newtheorem{lem}[df]{\sc Lemma}
\newtheorem{prp}[df]{\sc Proposition}
\newtheorem{ex}{\sc Example}


\title{\sc Growth of conjugacy classes of Schottky groups in higher rank
  symmetric spaces}
\author{\sc Gabriele Link}
\date{}
\maketitle
\begin{abstract} Let $\XX$ be a globally symmetric space of noncompact type,
  and $\Gamma\subset\is(\XX)$ a Schottky group of axial isometries. Then $M:=\XX/\Gamma$ is a locally
  symmetric Riemannian manifold of infinite volume. The goal of this note is
  to give an asymptotic estimate for the number of primitive closed
  geodesics in $M$ modulo free homotopy with  period less than $t$.
\end{abstract}
\vspace{0.2cm}

\section{Introduction}

Let $M$ be a complete Riemannian manifold of nonpositive sectional
curvature, and denote by $P(t)$ the number of primitive closed geodesics in
$M$ of period less than $t$ modulo free homotopy. If $M$ is compact with
volume entropy $h$, there are various results describing the asymptotic
behavior of this function $P(t)$: The most remarkable early result due to G.~A.~Margulis
(\cite{MargulisApplications},\cite{MR2035655}) states that if $M$ has pinched negative curvature, then
$$ \lim_{t\to\infty}P(t)\cdot h t\cdot  e^{-ht} =1\,.$$
Later, G.~Knieper (\cite{MR710022},\cite{MR1465601},\cite{MR1928523}) obtained a slightly weaker analogon
of this result for geometric rank one manifolds: He proved the existence of constants $a>1$ and $t_0>0$ \st
$$ \frac1{a\;t}\;e^{ht} \le P(t)\le \frac{a}{t}\; e^{ht}
$$ for $t>t_0$ (\cite[Theorem~5.6.2]{MR1928523}).  

For compact rank one symmetric spaces of noncompact type, Margulis' result has been improved by
giving error terms (\cite{MR461863}). However, for higher rank symmetric spaces not
even an analogon of Knieper's result is known. In this note we treat the case
where $M$ is a locally symmetric space of noncompact type 
with a Schottky group $\Gamma$ (in the sense of Y.~Benoist \cite{MR1437472}) as  deck transformation
group. Since $M$ has infinite volume, the exponential growth rate of $P(t)$
is  no longer governed by the volume entropy $h$ but instead by the {\hl
  critical exponent} of $\Gamma$ 
$$\delta(\Gamma):=\inf\{s>0\;|\, \sum_{\gamma\in\Gamma} e^{s d(x,\gamma
  y)}<\infty\}\,,$$ 
where $x$ and $y$ are arbitrary points in the Riemannian universal covering manifold
  of $M$.  
Our main result is the following \\[2mm]
{\bf Main Theorem}$\quad$ {\sl If $M$ is a locally symmetric space of rank $r\ge 1$ with deck transformation group $\Gamma$ as above, there 
exist constants $a>1$ and $t_0>0$ \st for any
$t>t_0$  
$$ \frac1{a\;t^r}\;e^{\delta(\Gamma)t}\le P(t)\le
\frac{a}{t}\;e^{\delta(\Gamma)t}\,.$$ }

We remark that for rank one symmetric spaces, the lower bound of the theorem holds in the more general context of convex cocompact deck transformation groups. This follows directly from Theorem~1 in \cite{Link2005Asymptotic}.

The paper is organized as follows: In section~2 we recall some
basic facts about symmetric spaces of noncompact type and decompositions of
semisimple Lie groups. Section~3  describes some important concepts concerning closed geodesics in a locally symmetric space. In section~4 we introduce the Schottky groups we will be concerned with and state some results about their limit set and the exponential growth rate of certain orbit points. 
Section~5 is devoted to the proof of the key step in our main theorem, namely to give an upper bound for the number of isometries corresponding to the same free homotopy class of closed geodesics in the quotient.
Finally, in section~6, we restate the main theorem and complete its proof.

\section{Preliminaries  on symmetric spaces}\label{Prelim}

The purpose of this section is to introduce some terminology and notation, and to summarize some basic results about symmetric spaces 
of noncompact type (see also \cite{MR1834454}, \cite{MR823981}, \cite{MR1441541}) which we shall need.

Let $\XX$ be a simply connected symmetric space of noncompact type
with base point $\xo\in\XX$,  
$G=\is^o(\XX)$ the connected component of the identity, and $K\subset G$ the
isotropy subgroup of $\xo$ in $G$. Then $\XX$ is a manifold of nonpositive curvature, $G$ a semisimple
Lie group with trivial center,  
$K$ a maximal compact subgroup of $G$,  and we may write $\XX=G/K$. 
If $\gL$ and $\kL$ denote the Lie
algebras of $G$ and $K$, then the geodesic symmetry of $\XX$ at $\xo$
determines a Cartan decomposition $\gL=\kL\oplus \pL$, where $\pL$ is
identified with the  tangent space $T_\xo\XX$ of $\XX$ at $\xo$. Let
$\aL\subset \pL$ be a maximal abelian subalgebra, and $\aL^+\subset \aL$ an
open Weyl chamber with closure $\overline{\aL^+}$. A {\hl flat} in $\XX$ is a totally geodesic submanifold of the form $ge^{\aL}\xo$, $g\in G$. The decomposition \Cart\ is
called the {\hl Cartan decomposition} of $G$. 
\begin{df}\label{thetavec}
For $x,y\in\XX$ the unique vector $H\in\overline{\aL^+}$ with the property 
$x=g\xo$ and $y=g e^H\xo$ for some $g\in G$
is called the {\hd Cartan vector}  of the ordered pair of points  $(x,y)\in\XX\times\XX$ and will be denoted $H(x,y)$. 
\end{df}
Notice that the length of the Cartan vector $H(x,y)$ is exactly the Riemannian distance between $x$ and $y$. If rank$(\XX):=\dim \aL=1$, then the Cartan vector of a pair of points 
is simply this number. 

Let $\Sigma$ be the set of restricted roots of the pair $(\gL,\aL)$, and
$\Sigma^+\subset\Sigma$ the set of positive roots determined by the
Weyl chamber $\aL^+$. We denote $\gL_\alpha$ the root space of
$\alpha\in\Sigma$, $\nL^+:=\sum_{\alpha\in\Sigma^+} \gL_\alpha$,
and $N^+$ the Lie group exponential of the nilpotent Lie algebra $\nL^+$. 
The decomposition $G=KAN^+$ is called the {\hl Iwasawa
  decomposition} associated to the Cartan decomposition \Cart. 
If $M$ denotes the centralizer of $\aL$ in $K$, the Iwasawa decomposition
induces a natural projection
\be \piKM\;:\quad G &\to & K/M\\
g=kan &\mapsto & kM\,. \ee 

Let $M^*$ be the normalizer of $\aL$ in $K$, and $W=M^*/M$ the {\hl Weyl group}
of the pair $(\gL,\aL)$.  We will denote $w_*\in W$ the unique element \st
$\Ad(m_{w_*})(-\aL^+)=\aL^+$ for any representative $m_{w_*}$ of $w_*$ in
$M_*$.

The geometric boundary $\rand$ of $\XX$ is defined as the set of equivalence
classes of asymptotic geodesic rays endowed with the cone
topology. This boundary is homeomorphic to the unit tangent space of
an arbitrary point in $\XX$ (compare \cite[chapter~II]{MR1377265}),
and $\ganz:=\XX\cup\rand$ is homeomorphic to a closed ball in $\RR^N$, where
$N:=\dim\XX$.  

Let $\aL_1\subset\aL$ be the set of unit vectors in $\aL\subset\pL\cong T_\xo
\XX$. For $x\in\XX$ and  $z\in\ganz\setminus\{x\}$ we denote $\sigma_{x,z}$
the unique unit speed geodesic ray emanating from $x$ containing $z$. The {\hl
  direction}  $H\in\overline{\aL^+_1}$ of a point $\xi\in\rand$ is defined as
the Cartan vector of the ordered pair $(\xo,\sigma_{\xo,\xi}(1))$,
i.e. $H=H(\xo,\sigma_{\xo,\xi}(1))$. If $k\in K$ is \st $\xi$ belongs to the
class of the geodesic ray $\sigma$ given by $\sigma(t):=k e^{Ht}\xo$, $t\ge
0$, we write $\xi=(k,H)$. 
Notice that $k$ is only determined up to right multiplication by an element in the centralizer of $H$ in $K$. 

If the rank of $\XX$ is greater than one, then the regular boundary $\regrand$ is defined as the set of
classes with Cartan projection $H\in\aL^+_1$. If $\rank(\XX)=1$, then
$\overline{\aL^+_1}$ is a point and we use the convention $\regrand:=\rand$. We will further need the 
continuous projection 
$$ \hspace{1cm} \ba{rcl}\pi^B\,:\qquad\ \regrand&\to& K/M\\
(k,H)&\mapsto & kM\,, \ea $$
which is a homeomorphism if and only if $\rank(\XX)=1$.

The isometry group of $\XX$ has a natural action by homeomorphisms on the
geometric boundary. If $g\in G$, $\xi=(k,H)\in\rand$ and $k'\in K$ is \st
$\pi^I(gk)=k'M$,   then $g\at (k,H) = (k', H)$ (see
\cite[Lemma~2.2]{Link2005Geometry}). In particular, the $G$-action preserves
the directions of boundary points, hence $G$ acts 
transitively on the geometric boundary if and only if $\rank(\XX)=1$.
However, the projection $\pi^B$ induces a transitive action
of $G$ by homeomorphisms on the {\hl Furstenberg boundary}
$K/M=\pi^B(\regrand)$. So if $\xi=(k,H) \in\regrand$, then
$g\pi^B(\xi)=\pi^B(g\xi)=k'M$. 

Moreover, for $\xi\in\rand$ we denote $\horinf(\xi)$ the set of points in the geometric boundary which can be joined to $\xi$ by a geodesic, i.e.
$$\horinf(\xi):=\{\eta\in\rand\;|\ \exists\ \mbox{geodesic}\ \sigma\  \st 
\sigma(-\infty)=\xi\,,\,\sigma(\infty)=\eta\}\,.$$ 
Notice that for rank one symmetric
spaces we have $\horinf(\xi)=\rand\setminus\{\xi\}$ for all $\xi\in\rand$. 
The {\hl Bruhat visibility set} of a point $\xi\in\regrand$
$$ \horF(\xi):=\pi^B\big(\horinf(\xi)\big)$$
will play an important role in the sequel. 
It is a dense and open submanifold of the Furstenberg boundary $K/M$ and
corresponds to a Bruhat cell of maximal dimension (compare
\cite[Section~2.3]{Link2005Geometry}). 

\section{Conjugacy classes and closed geodesics}\label{ConjClassClosed}

Let $M$ be a locally symmetric space of noncompact type 
with universal Riemannian covering manifold $\XX$, and $\Gamma\subset\is(\XX)$ the group of deck transformations of the covering projection $\XX\to M$. It is well-known that 
$\Gamma$ is a discrete and torsion free group isomorphic to the fundamental group of $M$. If $G=\is^o(\XX)$ is the connected component of the identity, we fix a Cartan decomposition \Cart\ and the associated Iwasawa decomposition $G=KAN^+$. Let $\xo\in\XX$ be the unique point stabilized by $K$. 

In order to describe closed geodesics in $M$, we will need to work with the following   kind of isometries of $\XX$:
\begin{df}  
An isometry $\gamma\neq\id$ 
of $\XX$ is called {\hd axial}, if there exists a constant
$l>0$ and a unit speed geodesic $\sigma\subset\XX$ \st
$\gamma(\sigma(t))=\sigma(t+l)$ for all $t\in\RR$.

We call $L(\gamma):=H(\sigma(0),\sigma(l))\in\overline{\aL^+}\setminus\{0\}$ the {\hd
    translation vector}, and $l(\gamma):=\Vert L(\gamma)\Vert>0$ the    {\hd
    translation length} of $\gamma$. 
The boundary point
$\gamma^+:=\sigma(\infty)$  is called the {\hd attractive fixed
point}, and $\gamma^-:=\sigma(-\infty)$ the {\hd repulsive fixed
point} of $\gamma$. 
We say that $\gamma$ is {\hd regular axial} if $\gamma^+\in\regrand$.
\end{df}
For an axial isometry $\gamma$ we define the set 
$\Ax(\gamma):=\{ x\in\XX\;|\, d(x,\gamma x)=l(\gamma)\}$ which consists of the union of parallel geodesics
translated by $\gamma$. Then 
$\overline{\Ax(\gamma)}\cap\rand$ is exactly the set of fixed points of
$\gamma$.

If $\gamma$ is regular axial and $x\in\Ax(\gamma)$, then there exists $g\in G$ \st
$x=g\xo$ and $\gamma x= g e^{L(\gamma)}\xo$. Moreover,
$\piF(\gamma^+)=\pi^I(g)$, $\piF(\gamma^-)=\pi^I(gw_*)$, and the set of fixed
points in the Furstenberg boundary $\Fix^B(\gamma)$ is exactly the finite set 
$$\Fix^B(\gamma)= \{ \pi^I(gw)\;|\, w\in W\}\,.$$

If $\gamma\in\Gamma$ is axial and $\sigma$, $\sigma'\subset\Ax(\gamma)\subset\XX$ are geodesics translated by $\gamma$, then $\sigma$ and $\sigma'$ project to  freely homotopic closed geodesics of the same period $\le l(\gamma)$ in the quotient $M$. Since $\Ax(\gamma)$ consists of an uncountable union of parallel geodesics translated by $\gamma$, there are uncountably many closed geodesics in each free homotopy class.  In order to describe these free homotopy classes, we will make use of the following
\begin{df} $\gamma, \gamma'\in\Gamma$ are said to be {\hd
equivalent} if and only if there exist $n,m\in \ZZ$ and
$\varphi\in\Gamma$ \st $(\gamma')^m=\varphi \gamma^n\varphi^{-1}$.
An element $\gamma_0\in\Gamma$ is called {\hd primitive} if it
cannot be written as a proper power $\gamma_0=\varphi^n$, where
$\varphi\in\Gamma$ and $n\ge 2$. 
\end{df} 
Each equivalence class
can be represented as
$$[\gamma]=\{\varphi \gamma_0^k\varphi^{-1}\;|\,
\gamma_0\in\Gamma,\;\gamma_0\
\mbox{primitive},\;k\in\ZZ,\;\varphi\in\Gamma\}\,.$$ 
It is easy to see that
the set of equivalence classes of axial elements in $\Gamma$ is in
one to one correspondence with the set of geometrically distinct
closed geodesics modulo free homotopy.
If $\gamma\in\Gamma$ is axial, we put 
$$l([\gamma]):=\min \{l(\varphi)\;|\, \varphi\in [\gamma]\}$$
and notice that if $\gamma_0$ is a primitive isometry representing
$[\gamma]$, then $l([\gamma])= l(\gamma_0)$.  
Moreover,  
$$P(t):=\#\{[\gamma]\;|\, \gamma\in\Gamma\ \mbox{axial},\;
l([\gamma]) < t\}$$ 
counts the number of geometrically distinct closed geodesics of period less
than $t$ modulo free homotopy.

We will see that $P(t)$ is intimately related to the number
$$N_\Gamma(R):=\#\{\gamma\in\Gamma\;|\,d(\xo,\gamo)<R\}\,,$$
which allows the following alternative characterization of the critical exponent of $\Gamma$:
$$\delta(\Gamma)=\limsup_{R\to\infty}\big( \frac1{R}\log N_\Gamma(R)\big)\,.$$
In the sequel, we will further need the following 
\begin{df}\label{VarLimSets}
If  $\Gamma\subset\is(\XX)$ is a discrete group, its {\hl limit set} $\Lim$ is defined
by  $\Lim:=\Gamma\at\xo\cap \rand$. 
We put $\Flim:=\piF(\Lim\cap\regrand)$,
and $P_\Gamma:=\{ H\in \overline{\aL_1^+}\;|\, \exists\; k\in K\ \,
\mbox{such that}\ (k,H)\in\Lim\}$.
\end{df}

\section{Schottky groups}\label{Schottky}

We next introduce and recall some properties of the Schottky groups we will be
concerned with in the sequel. As in the previous section, $\XX$ is a globally symmetric
space of noncompact type, \Cart\  a Cartan decomposition of $G=\is^o(\XX)$,
and $\xo\in\XX$ the unique point stabilized by $K$. 

Let $\gax_1, \gax_2,\ldots, \gax_l$ be regular axial
isometries with the following property: If we denote $\xi_{2m}:=\gax_m^+$ and
$\xi_{2m-1}:=\gax_m^-$, $1\le m\le l$, then
$$  \pi^B(\xi_i)  \in  \bigcap_{\begin{smallmatrix}n=1\\
n\neq i\end{smallmatrix}}^{2l} \horF(\xi_n)\qquad\forall\, i\in\{1,2,\ldots,
2l\}\,.$$
Up to replacing $\gax_1$ and $\gax_2$ by approximate elements, we may assume
that the group generated by $\gax_1, \gax_2,\ldots, \gax_l$ is Zariski
dense in $G$ (compare \cite[Section~4.2]{MR2147895}). Furthermore, by
Lemma~4.5.15 in \cite{MR1441541}  there exist
neighborhoods $U'_n\subset K/M$ of $\piF(\xi_n)$, $1\le n\le 2l$, such that for all $i\in \{1,2,\ldots, 2l\}$ and every point $\eta\in\regrand$ 
with $\piF(\eta)\in U'_i $ we have
\begin{equation}\label{oppositioncondition}
\bigcup_{\begin{smallmatrix}n=1\\
n\neq i\end{smallmatrix}}^{2l} U'_n \subset \horF(\eta)\,.
\end{equation}
 After possibly replacing $\gax_1,\gax_2,\ldots,\gax_l$ by
sufficiently large powers, we may assume that there exist sets 
$(U_n,W_n)_{1\le n\le 2l}$ in $K/M$, $U_n\subset U_n'$,  \st conditions (i), (ii), (iii) before
Proposition~4.4 in \cite{MR2147895} are satisfied with  $b_i^+=U_{2i}$,
$b_i^-=U_{2i-1}$, $B_i^+=W_{2i}$, $B_i^-=W_{2i-1}$, $1\le i \le l$, for some
$\eps\in (0,1)$. In particular, 
$$\bigcap_{n=1}^{2l} W_n\neq \emptyset\,, \qquad U_{2i}\subset \bigcap_{\begin{smallmatrix}n=1\\
n\neq 2i-1\end{smallmatrix}}^{2l} W_n\,, \qquad U_{2i-1}\subset \bigcap_{\begin{smallmatrix}n=1\\
n\neq 2i\end{smallmatrix}}^{2l} W_n \,,\qquad \mbox{and}$$
$$ \gax_i(W_{2i})\subset U_{2i}\,,\quad \gax_i^{-1}(W_{2i-1})\subset U_{2i-1}
\qquad\mbox{for all}\quad   i\in \{1,2,\ldots,l\}\,.$$
Hence the group $\Gamma:=\langle \gax_1, \gax_2,\ldots, \gax_l\rangle$
is a free and discrete Zariski dense subgroup of $G$. Moreover, we have the following well-known 
\begin{prp}(\cite{MR1437472}, \cite[Section~4.5]{Link2005Geometry})\label{LimSetDescription}\\
The  limit set of $\Gamma$ is contained in the regular boundary and splits as a product 
$\Lim\cong\Flim\times P_\Gamma$. Moreover, $\Flim$ is a minimal closed set for the
action of $\Gamma$, $\Flim\subset \bigcup_{i=1}^{2l} \overline{U_i}$ and
$P_\Gamma\subset\aL^+_1$ is a closed convex cone. Every element $\gamma\in\Gamma$ is regular axial and satisfies $L(\gamma)/l(\gamma)\in P_\Gamma$. 
\end{prp} 
 
In particular, Theorem~5.1 in \cite{MR2147895} applies to $\Gamma$,  hence there exist constants $b>1$ and $R_0>0$ \st for all $R>R_0$
\begin{equation}\label{quintresult}
 \frac1{b} e^{\delta(\Gamma)R}\le N_\Gamma(R)\le b e^{\delta(\Gamma)R}\,.
\end{equation} 

As one of the main ingredients in the  proof of the lower bound will serve Proposition~\ref{countwithsets}, a stronger version of
equation~(\ref{quintresult}). The idea of proof is originally  due to T.~Roblin
(\cite{MR1881574}), but, due to
the necessity of dealing both with the geometric boundary and the Furstenberg
boundary, is more technical in our situation. For this reason the following
definitions  will be convenient. 

If $C\subset K/M$, $z\in\ganz$, we put
\begin{equation}\label{AngleDef}
\angle_o(z,C):=\inf\{\angle_o(z, \eta)\;|\, \eta\in\regrand\ \, \mbox{such
  that}\  \piF(\eta)\in C\}\,,
\end{equation}
for $A,B\subset K/M$ we denote 
$$N_\Gamma(R;A,B):=\#\{\gamma\in\Gamma\;|\, d(\xo,\gamo)<R,\;
\angle_\xo(\gamo,A)=0,\; \angle_\xo(\gamma^{-1}\xo,B)=0\}\,.$$
If $V\subset K/M$ and $\eps>0$ we define a subset of the regular boundary by  
\begin{eqnarray}\label{NeighborhoodDef}
 H^\eps_\Gamma (V)&:= &\{\eta\in\regrand\;|\, \exists\ \xi=(k,H)\in\regrand\ \,
  \mbox{with}\ \, kM\in V\ \mbox{and}\ H\in P_\Gamma\nonumber\\
&&\hspace{6.3cm} \mbox{such
  that} \ \ \angle_o(\eta, \xi) < \eps\}\,.
\end{eqnarray}
If $C\subset K/M$ is an open set containing the closure $\overline V$ of $V$,
we further put
\begin{equation}\label{epsinradius}
\eps_\Gamma(V,C):=\sup\{\eps >0\;|\, \piF(V_\eps(\Gamma))\subset C\}\,.
\end{equation} 
Moreover, the following easy lemma will be necessary in the proof of Proposition~\ref{countwithsets}.
\begin{lem}\label{impossibleposition}
Let $V\subset C\subset K/M$ be as above and put $\eps:=\eps_\Gamma(\overline{V},C)$. Then
for all $z\in\ganz$ with $\angle_\xo(z,C)>0$ and $\angle
(H(\xo,z),P_\Gamma)<\eps/3$ we have
$$ \angle_\xo(z,\xi)>\frac{\eps}3\quad\quad\forall\, \xi\in\regrand\ \mbox{
  with}\quad \piF(\xi)\in\overline{V}\,.$$
\end{lem}
\prf\ Suppose $\xi\in\regrand$ satisfies $\piF(\xi)\in\overline{V}$ and
$\angle_\xo(z,\xi)\le\eps/3$. If $H\in\aL^+_1$ denotes the
direction of $\xi$, then $\angle (H(\xo,z),H)\le \angle_\xo(z,\xi)\le
\eps/3$. Hence $\angle(H(\xo,z),P_\Gamma)<\eps/3$ implies
$$\angle (H,P_\Gamma)\le \angle (H,H(\xo,z))+\angle
(H(\xo,z),P_\Gamma)<\frac23\eps\,.$$
In particular, if $\xi=(k,H)$ and $H'\in P_\Gamma$ \st
$\angle(H,P_\Gamma)=\angle (H,H')$, then $\eta:=(k,H')\in\regrand$ satisfies 
 $\angle_\xo (\xi,\eta)=\angle (H,H')<2\eps/3$. We conclude that
$$ \angle_\xo(z,\eta)\le \angle_\xo(z,\xi)+\angle_\xo(\xi,\eta)<\eps\,.$$
Moreover, since $\piF(\eta)=\piF(\xi)\in\overline{V}$ and $H'\in P_\Gamma$ we
have $\sigma_{\xo,z}(\infty)\in H_\Gamma^\eps(\overline{V})$, hence by choice
of $\eps=\eps_\Gamma(\overline{V},C)$ $\ \piF\big(
\sigma_{\xo,z}(\infty)\big)\in C$, in contradiction to
$\angle_\xo(\sigma_{\xo,z}(\infty),C)=\angle_\xo(z,C)>0$.\qed\\
\begin{prp}\label{countwithsets} 
If $A,B \subset K/M$ are open sets with $\piF(\gamma^+)\in A$ and $\piF(\gamma^-)\in B$ for some element $\gamma\in\Gamma$, then 
there exist constants $a>1$ and $R_0>0$ \st for all $R>R_0$ 

$$ \frac1{a}\; e^{\delta(\Gamma)R}\le N_\Gamma(R;A,B)\le a\; e^{\delta(\Gamma)R}\,.$$ 
\end{prp}
\prf\ Let $U, V \subset K/M$ be  open sets with $\piF(\gamma^+)\in U$, $\piF(\gamma^-)\in V$, $\overline{U}\subset A$ and $\overline{V}\subset B$. 

We first prove the claim for $N_\Gamma(R;A):=\#\{\gamma\in\Gamma\;|\, d(\xo,\gamo)<R,\; \angle_\xo(\gamo,A)=0\}$. 
Since $\Flim$ is compact, and $\Flim=\piF(\overline{\Gamma\at\gamma^+})$
by Proposition~\ref{LimSetDescription}, there exist
$g_1,g_2,\ldots,g_m\in\Gamma$ \st
$\Flim\subseteq\bigcup_{i=1}^m g_i U$. Let 
$Y\subseteq\bigcup_{i=1}^m g_i U$ be an open set which
contains $\Flim$. Then $\angle_o(\gamo,Y)=0\,$ for all but finitely
many $\gamma\in\Gamma$ by the definition of the limit set and (\ref{AngleDef}), hence
$N_\Gamma(R; Y)\ge N_\Gamma(R)-M$ for some constant $M\in\NN$.

Fix $g\in\{g_1,g_2,\ldots,g_m\}$ and suppose $\angle_\xo(\gamo,gU)=0$ and $\angle_\xo(g^{-1}\gamo,A)>0$ for infinitely many
$\gamma\in\Gamma$. Let $(\gamma_j)\subset\Gamma$ be a sequence
with this property. Passing to a subsequence if necessary, we may
assume that $\gamma_j\xo$ converges to a point $\eta\in \Lim\subset\regrand$, and we
have $\piF(\eta)\in g\overline U$. Let $\eps:=\eps_\Gamma(\overline{U},A)$ 
and $N_0\in\NN$ \st for all $j\ge N_0$ we have 
\begin{equation}\label{convergence}
\angle_\xo(g^{-1}\gamma_j\xo,g^{-1}\eta)
=\angle_{g\xo}(\gamma_j\xo,\eta)<\eps/4\,. 
\end{equation}
On the other hand, by choice of $(\gamma_j)$ and since
$\piF(g^{-1}\eta)\in\overline{U}$, Lemma~\ref{impossibleposition} implies\\   
$\angle_\xo (g^{-1}\gamma_j\xo, g^{-1}\eta)>\eps/3$ for all $j$ sufficiently
large, in contradiction to (\ref{convergence}).

We conclude that
$$c(g):=\#
\{\gamma\in\Gamma\;|\, \angle_{\xo}(\gamo,g U)=0\ \an\
\angle_o(g^{-1}\gamo,A)>0\}$$ 
is finite and  
$N_\Gamma(R;g U)\le
N_\Gamma(R+d(\xo,g\xo);A)+c(g)\,.$ With $d:=\max_{1\le i\le
m} d(\xo,g_i\xo)$ and $c:=\sum_{i=1}^m c(g_i)$  we obtain 
\be
N_\Gamma(R-d;A)-M & \le & N_\Gamma(R-d)-M  \le
N_\Gamma(R-d;Y)\le\sum_{i=1}^m N_\Gamma(R-d;g_i U)\\
&\le & \sum_{i=1}^m \big( N_\Gamma(R-d+d(\xo,g_i\xo);A)+c(g_i)\big)\le mN_\Gamma(R;A)+c\,,
\ee 
which proves the first assertion. 

Now again by compactness of $\Flim$ and the fact that $\Flim=\piF(\overline{\Gamma\at\gamma^+})$, there exist
$g_1,g_2,\ldots,g_n\in\Gamma$ \st
$\Flim\subseteq\bigcup_{i=1}^n g_i V$. By the same argument as above the numbers 
$$ c(g) :=\#\{\gamma\in\Gamma\;|\, \angle_{\xo}(\gamma^{-1}\xo,g V)=0\ \an\
\angle_\xo(g^{-1}\gamma^{-1}\xo,B)>0\}$$ are finite for all $g\in\{g_1, g_2,\ldots,g_n\}$, and
$N_\Gamma(R;A,g V)\le
N_\Gamma(R+d(\xo,g\xo);A,B)+c(g)\,.$ We put $d:=\max_{1\le i\le
n} d(\xo,g\xo)$ and $c:=\sum_{i=1}^n c(g_i)$ and conclude
$$N_\Gamma(R-d;A,B)\le N_\Gamma(R-d;A)\le
\sum_{i=1}^n N_\Gamma(R-d;A,g_i V)\le n
N_\Gamma(R;A,B)+c\,,$$ which yields the assertion.\qed\\

\section{The key step}\label{lowbound}

The most difficult task when trying to obtain a lower bound for $P(t)$ is to
give an upper bound for the number of elements in the same  conjugacy
class. The main problem compared to the situation in \cite{MR1465601} or \cite{Link2005Asymptotic} arises from the fact that in the higher rank setting, a cyclic group generated by an axial isometry does not act cocompactly on its invariant set. In the case of Schottky groups as in Section~\ref{Schottky}, however, we are able to show that for every isometry $\gamma\in\Gamma$, the number of orbit points in $\XX$ close to such an unbounded fundamental set 
for the action of $\langle \gamma\rangle $ on $\Ax(\gamma)$ is bounded from above by a function of $l(\gamma)$. This will finally allow to obtain the desired estimate. 

We start with a couple of preliminary lemmata. For the remainder of this section,
$\Gamma=\langle \gax_1,\gax_2,\ldots,\gax_l\rangle\subset\is(\XX)$ will be a
Schottky group acting on a globally symmetric space $\XX$ of rank $r>1$ as in the
previous section.
Recall Definition~\ref{VarLimSets} 
and put 
$$\alpha_\Gamma:=\sup_{H,H'\in P_\Gamma} \angle (H,H')\,.$$ 
Since $P_\Gamma$ is a closed convex cone strictly included in $\aL_1^+$, there
exists $\delta\ge 0$ such that for any $H\in\aL$ 
\begin{equation}\label{deltadef}
\angle (H,P_\Gamma):=\inf_{H'\in P_\Gamma} \angle (H,H')\le
\delta\quad\mbox{implies}\ H\in\aL^+\,.
\end{equation} 
From the fact that $\max\{ \angle (H,H')\;|\, H,H'\in\overline{\aL^+}\}\le
\pi/2$ (see i.e. \cite[Theorem~X.3.6 (iii)]{MR1834454}) we conclude that $\alpha_\Gamma +\delta$ is strictly smaller than
$\pi/2$.

Given $\gamma\in\Gamma$, we denote $x_\gamma$ the orthogonal projection of
$\xo$ to $\Ax(\gamma)$. Let $g\in G$ \st $x_\gamma=g\xo$ and $\gamma
x_\gamma=g e^{L(\gamma)}\xo$, and put
\begin{eqnarray*}
F(\gamma) &:=&\{g  e^{H}\xo\;|\, H\in\aL\ \, \st\, \langle
H,\frac{L(\gamma)}{l(\gamma)}\rangle\in [0,l(\gamma)]\}\,,\\
C_\Gamma^\delta(\gamma)&:=& \{g  e^{H }\xo\;|\,  H\in\aL\ \, \st\, \angle (H, P_\Gamma)\le
\delta \}\subset g e^{\aL^+}\xo \,.
\end{eqnarray*}
Notice that $F(\gamma)$, $C_\Gamma^\delta(\gamma)\subset\Ax(\gamma)$ are closed sets,
and $\langle \gamma\rangle\cdot  F(\gamma)=\Ax(\gamma)$. We further have the following
\begin{lem}\label{volestimate}
There exists a constant $a=a(\Gamma,\delta,r)>0$ \st for all $\gamma\in\Gamma$ 
$$ {\rm vol}\big(F(\gamma)\cap C_\Gamma^\delta(\gamma)\big)\le a\at l(\gamma)^r\,.$$
Furthermore, $F(\gamma)\cap C_\Gamma^\delta(\gamma)$ is compact for any $\gamma\in\Gamma$. 
\end{lem}
\prf\ For $\gamma\in\Gamma$ we put 
$$\aL(\gamma):=\{ H\in\aL^+\;|\, \langle H, L(\gamma)\rangle\le l(\gamma)^2,
\, \angle (H, P_\Gamma)\le \delta\}\,.$$
Then $\vol\big(F(\gamma)\cap C_\Gamma^\delta(\gamma)\big)=\int_{\aL(\gamma)}dH$.

We substitute $\hat H=H/\Vert H\Vert$, $t=\Vert H\Vert$, and remark that $\angle  (H,
P_\Gamma)\le\delta$ and $L(\gamma)/l(\gamma)\in P_\Gamma$ imply 
$$\angle (H,L(\gamma))\le \alpha_\Gamma+\delta<\frac{\pi}2\,,\quad\mbox{hence}\quad \ \langle \hat H,L(\gamma)\rangle\ge l(\gamma)\at \cos (\alpha_\Gamma+\delta)>0 \,.$$
In particular, $H\in \aL(\gamma)$ satisfies the condition $\Vert H\Vert \at\cos
(\alpha_\Gamma+\delta)\le l(\gamma)$. We summarize
$$  \int_{\aL(\gamma)}dH\le
\int_0^{l(\gamma)/\cos(\alpha_\Gamma+\delta)}\underbrace{\big(\int_{\angle
    (\hat H,  P_\Gamma)\le\delta} d\hat
H\big)}_{=:v_\Gamma(\delta)} t^{r-1}dt=\frac{v_\Gamma(\delta)}{r\at (\cos(\alpha_\Gamma+\delta))^r}\cdot
l(\gamma)^r\,,$$ hence the assertion follows with $a=v_\Gamma(\delta)/(r\at
\cos^r(\alpha_\Gamma+\delta))$. 

Moreover, $F(\gamma)\cap C_\Gamma^\delta(\gamma)$ is compact since $\aL(\gamma)$ is
closed and the norm of every element in $\aL(\gamma)$ is bounded by
$l(\gamma)/\cos(\alpha_\Gamma+\delta)<\infty$. \qed\\

Using the notation from the previous section we recall that
$K_\Gamma\subset\bigcup_{i=1}^{2l}\overline{U_i}$, and $\overline{U_i}\subset U_i'$ for
all $1\le i\le 2l$. We put 
\be \Gamma'&:= &\{\gamma\in \Gamma\;|\, \piF(\gamma^-)\in
U_1\quad \mbox{and}\ \,\piF(\gamma^+)\in U_2\}\,,\ \; \mbox{and}\\
K'&:= &\{kM\in K/M\;|\, \exists\; \gamma\in\Gamma' \quad \st\ \, kM\in\Fix^B(\gamma)\setminus\{\piF(\gamma^-),\piF(\gamma^+)\}\}\,.\ee

$K'$ describes the set of fixed points in the Furstenberg boundary of elements
in $\Gamma'$ which do not correspond to attractive or repulsive fixed
points. Since for all $\gamma\in\Gamma$ we have $\horF(\gamma^-)\cap
\Fix^B(\gamma)=\{\piF(\gamma^+)\}$ and $\horF(\gamma^+)\cap
\Fix^B(\gamma)=\{\piF(\gamma^-)\}$, condition~(\ref{oppositioncondition}) implies that
$\ K'\subset K/M\setminus \bigcup_{i=1}^{2l}U_i'\subseteq K/M\setminus K_\Gamma$.

Moreover, by $\overline{U_i}\subseteq U_i'$ for $1\le i\le 2l$, the number
$$\min_{1\le i\le 2l} \eps_\Gamma(\overline{U_i}, U_i')$$ is positive. We fix $\eps\in
(0,\delta)$ strictly smaller than this minimum and put 
$$\Lim^\eps:=\{\eta\in\regrand\;|\, \exists\; \xi\in\Lim\ \st
\angle_\xo(\xi,\eta)<\eps\}\subseteq H_\Gamma^\eps (K_\Gamma)\subseteq
\bigcup_{i=1}^{2l} H_\Gamma^\eps (\overline{U_i})\,.$$
Then by choice of $\eps$ we have $\piF(\Lim^\eps)\cap \overline{K'}=\emptyset$,  hence there can be only finitely many orbit points close to $K'$. More precisely, we have 
\begin{lem}
If $\eps\in (0,\delta)$ is as above, then $\ N_\eps:=\#\{\gamma\in\Gamma\;|\, \angle_\xo (\gamma\xo,K')<\eps/4\}<\infty\,.$
\end{lem}
\prf\ Suppose $N_\eps$ is infinite. Then there exists a sequence
$(\gamma_j)\subset\Gamma$ \st $\angle_\xo(\gamma_j\xo,K')<\eps/4$. 
Let $(\eta_j)\subset\regrand$, $\piF(\eta_j)\in K'$, \st
$\angle_\xo(\gamma_j\xo,\eta_j)<\eps/2$ for all $j\in\NN$. Passing to subsequences if necessary,
we may assume that $\gamma_j\xo\to \xi\in\Lim\cap\regrand$ and
$\eta_j\to \eta\in\rand$. 
Hence for $j$ sufficiently large we have
$\angle_\xo(\gamma_j\xo,\xi)<\eps/4$ and $\angle_\xo(\eta_j,\eta)<\eps/4$, and
we conclude
$$ \angle_\xo(\xi,\eta)\le
\angle_\xo(\xi,\gamma_j\xo)+\angle_\xo(\gamma_j\xo,\eta_j)+\angle_\xo(\eta_j,\eta)<\frac{\eps}4+\frac{\eps}2+\frac{\eps}4=\eps\,.$$
If $H_\xi\in P_\Gamma$, $H_\eta\in\overline{\aL_1^+}$ are the directions of
$\xi$ and $\eta$, this implies $\angle (H_\xi,H_\eta)\le
\angle_\xo(\xi,\eta)<\eps<\delta$, hence by choice of $\delta$ we have 
$\eta\in\regrand$. From $\piF(\eta_j)\in K'$ and the definition of
$\Lim^\eps$ we conclude $\piF(\eta)\in \overline{K'}\cap \piF(\Lim^\eps)$, 
a contradiction to $\piF(\Lim^\eps)\cap \overline{K'}=\emptyset$.\qed\\

We next put 
$$c_0:=\sup \{ d(\xo, \Ax(g)) \;|\, g\in\is(\XX)\ \mbox{regular axial with}\ \piF(g^-)\in U_1,\, \piF(g^+)\in U_2
\}\,$$
and recall Definition~\ref{thetavec} for the Cartan vector of a pair of points in $\XX$. 
Let 
\begin{equation}\label{rhodef}
\rho:=\frac14 \min_{\gamma\in\Gamma}\big(\inf_{x\in B_\xo(c_0)}
d(x,\gamma x)\big)\,
\end{equation}
and choose $R=R(\eps,c_0,\rho)>0$ \st the following three
conditions hold:
\begin{enumerate}
\item[(i)] $\forall\, x,y\in B_\xo(c_0)\quad\forall\;\gamma\in\Gamma \;: \quad
  d(y,\gamma x)>R-c_0-\rho\  \Longrightarrow\ \angle (H(y,\gamma x),P_\Gamma) <\eps/8$,
\item[(ii)] $\forall\, x,p,z\in\XX\;: \quad d(x,p)>R-c_0\,,\ d(p,z)<\rho+c_0\
  \Longrightarrow\ \angle_x(p,z)<\eps/8$,
\item[(iii)] $\forall\, x,p\in\XX\;: \quad d(\xo,p)>R\,,\ d(\xo,x)<c_0\
  \Longrightarrow\ \angle_\xo(p,\sigma_{x,p}(\infty))<\eps/8$.
\end{enumerate}
Our first observation is the following 
\begin{lem}\label{CartVektorEstimate}
Let $x\in B_\xo(c_0)$, $\gax\in\Gamma'$ and $x_\gax\in \Ax(\gax)$ the
orthogonal projection of $\xo$ to $\Ax(\gax)$. Then for any $p\in \Ax(\gax)$ with
$d(\xo,p)>R$ and $B_p(\rho)\cap\Gamma\at x\neq \emptyset$ we have
$\angle (H(x_\gax, p), P_\Gamma)<\eps/4$.
\end{lem}
\prf\ Since $B_p(\rho)\cap\Gamma\at x\neq \emptyset$, there exists
$\gamma\in\Gamma$ \st $d(p, \gamma x)<\rho$. By choice of
$c_0>0$ we have $d(\xo, x_\gax)<c_0$, hence the triangle inequality implies
$$ d(x_\gax,\gamma x)\ge d(\xo,p)-d(\xo,x_\gax)-d(p,\gamma x) >R-c_0-\rho\,.$$
Condition~(i) above then gives $\angle (H(x_\gax,\gamma x),P_\Gamma)<\eps/8$,
and  $d(x_\gax,p)>R-c_0$ and $d(p,\gamma x)<\rho$ imply
$\angle_{x_\gax}(p,\gamma x)<\eps/8$ by condition~(ii).  We estimate 
$$\angle(H(x_\gax,p),H(x_\gax,\gamma x))\le \angle_{x_\gax}(p,\gamma
x)<\eps/8\,$$
and conclude $\angle (H(x_\gax,p),P_\Gamma)<\eps/4$. \qed\\

We will now combine these lemmata in order to obtain the desired estimate.
\begin{prp}
Let $\eps\in (0,\delta)$, $c_0>0$ and $\rho>0$  as above. 
Then there exist constants $b>0$ and $t_0>0$ 
\st for any $x\in \overline{B_\xo(c_0)}$,  every
primitive element $\gax\in\Gamma'$ and all $\,t>t_0$ 
\begin{eqnarray*}
\#\{\gamma=\varphi \gax^k\varphi^{-1}\;|\,\varphi\in\Gamma,\;k\in\ZZ,\;
\piF(\varphi\gax^-)\in U_1,\,\piF(\varphi\gax^+)\in U_2,\,\ && \\
\varphi\Ax(\gax)\cap B_x(\rho)\cap B_\xo(c_0)\ne\emptyset,\,\  l(\gamma)\le t \} & \le &  b \cdot t^r\,.
\end{eqnarray*}
\end{prp}
\prf\ We first remark that if $\gamma=\varphi \gax^k\varphi^{-1}$ with
$l(\gamma)\le t$, then $t\ge l(\gax^k)=|k|\cdot l(\gax)$, hence $|k|\le t/l(h)$.

Now if $k\in\ZZ\setminus\{0\}$ is fixed, then
$\varphi\gax^k\varphi^{-1}\ne \beta\gax^k\beta^{-1}$
implies that $\beta^{-1}\varphi$ does not belong to the
centralizer $Z_\Gamma(\gax)$ of $\gax$ in $\Gamma$, in particular
$\varphi\Ax(\gax)\neq\beta\Ax(\gax)$. If
$$\Gamma(\gax):=\{\gamma\in\Gamma\;|\, \piF(\gamma\gax^-)\in U_1,\,\  \piF(\gamma\gax^+)\in U_2,\,\ \gamma\Ax(\gax)\cap
B_x(\rho)\cap B_\xo(c_0)\ne\emptyset\}\,,$$ then the number of different conjugates $\gamma=\varphi \gax^k\varphi^{-1}$ with $\varphi\in\Gamma$, $\piF(\varphi\gax^-)\in U_1,\,\  \piF(\varphi\gax^+)\in U_2$ and $\varphi\Ax(\gax)\cap
B_x(\rho)\cap B_\xo(c_0)\ne\emptyset$ equals the cardinality of $\Gamma(\gax)$ modulo the
centralizer $Z_\Gamma(\gax)$ of $\gax$ in $\Gamma$. 

Denote  $x_\gax\in\Ax(\gax)$ the projection of $\xo$ to $\Ax(\gax)$, let
$g_\gax\in G$ \st $x_\gax=g_\gax \xo$ and $\gax x_h=g_\gax 
e^{L(\gax)}\xo$, and  
$$F(\gax)=\{ g_\gax e^H\xo\;|\, H\in\aL\ \, \st\, \langle
H, L(\gax)\rangle\in [0,l(\gax)^2]\} \subset \Ax(\gax)$$ 
as defined before Lemma~\ref{volestimate}.
If $\varphi,\beta\in \Gamma(\gax)$, $\beta^{-1}\varphi\notin Z_\Gamma(\gax)=\langle \gax\rangle$, there 
exist $p,q\in F(\gax)$ and $n,m\in\ZZ$ \st
$\varphi\gax^n p$, $\beta\gax^m q\in B_x(\rho)\cap B_\xo(c_0)$. Further\-more
$g:=\varphi\gax^n\ne \beta\gax^m=:f$ implies $d(p,q)\ge 2\rho$ since 
\be 2 \rho&\ge & d(gp,fq)\ge
d(gp,fp)-d(fp,fq)\\
&=& d(gf^{-1}gp,gp)-d(p,q)\stackrel{(\ref{rhodef})}{\ge} 4 \rho-d(p,q)\,.\ee 
So if we assign to each element $\varphi\in\Gamma(\gax)/\langle\gax\rangle$ a unique point
$p=p(\varphi)\in F(\gax)$ as above, then the balls $B_{p(\varphi)}(\rho)$,
$\varphi\in\Gamma(\gax)/\langle \gax\rangle$, are pairwise disjoint. 

For $R>0$ as before Lemma~\ref{CartVektorEstimate} we are going to bound the cardinality of
$$ Y_R(\gax):=\{ p(\varphi)\;|\, \varphi\in \Gamma(\gax)/\langle\gax\rangle,\,
d(\xo,p(\varphi))>R\}\subset F(\gax)\,.$$
If $p\in Y_R(\gax)$, then $p\in F(\gax)$ and we may write $p=g_\gax e^H\xo$
with $H\in \aL$, $\langle H,L(\gax)\rangle \in [0,l(\gax)^2]$. Furthermore, $H=\Ad(w)H(x_\gax, p)$ for
some $w\in W$, and we have $w\ne w_*$ by the fact that $\langle
-H',L(\gax)\rangle <0$ for all $H'\in \aL^+\setminus\{0\}$. Hence the 
following two cases may occur:
\begin{itemize}
\item[1.] Case: $\quad H=H(x_\gax,p)\in \aL^+$ or, equivalently, $w=\id$:\\
By Lemma~\ref{CartVektorEstimate} we have $\angle (H,P_\Gamma)<\eps/4<\delta$, hence
$p\in C_\Gamma^{\delta}(\gax)$. Since the balls $B_{p(\varphi)}(\rho)$,
$\varphi\in\Gamma(\gax)$, are pairwise disjoint, and $ \vol(B_p(\rho)\cap
\Ax(\gax))=\omega_r\cdot\rho^{r}$ for any $p\in\Ax(\gax)\,,$ there are at most 
$$\frac{\vol\big(F(\gax)\cap C_\Gamma^\delta(\gax)\big)}{\omega_r\cdot
  \rho^r}\le \frac{a\at l(\gax)^r}{\omega_r\cdot\rho^{r}}$$ different
points $p(\varphi)\in Y_R(\gax)$ with $p(\varphi)=g_\gax
  e^{H(x_\gax,p(\varphi))}\xo$. 
\item[2.] Case: $\quad H=\Ad(w)H(x_\gax,p)$ for some $w\in W\setminus\{\id,w_*\}$:\\
We have $p=g_\gax we^{H(x_\gax,p)}\xo$ and put
$\xi:=\sigma_{x_\gax,p}(\infty)\in\regrand$. Then
$\piF(\xi)\in\Fix^B(\gax)\setminus\{\piF(\gax^+),\piF(\gax^-)\}$, hence
$\gax\in\Gamma'$ implies
$\piF(\xi)\in K'$. Let $\gamma\in\Gamma$ \st $\gamma x\in B_p(\rho)$. By condition~(ii) above,
$d(\xo,p)>R$ and $d(p,\gamma\xo)\le d(p,\gamma x)+d(\gamma x,\gamma \xo)<\rho+c_0$ imply $\angle_\xo
(p,\gamma\xo)<\eps/8$. In particular, 
$$\angle_\xo(\gamma\xo, K')\le \angle_\xo(\gamma\xo,\xi)\le
\angle_\xo(\gamma\xo,p)+\angle_\xo(p,\xi)\stackrel{\mbox{{\scriptsize (iii)}}}{<}\frac{\eps}8+\frac{\eps}8=\frac{\eps}4\,.$$
This shows that there are at most $N_\eps$ points $p(\varphi)\in Y_R(\gax)$
of the form\\
$p(\varphi)=g_\gax w  e^{H(x_\gax,p(\varphi))}\xo$ with $w\ne\id$.  
\end{itemize}
We summarize $\ \# Y_R(\gax)\le a\cdot
l(\gax)^r/(\omega_r\cdot\rho^{r})+N_\eps\,,$ hence
$$ \#\ \Gamma(\gax)/\langle\gax\rangle\le
N_\Gamma(R+\rho)+\frac{a}{\omega_r\cdot\rho^{r}}\cdot l(\gax)^r+N_\eps\,.$$
Since $ N_\Gamma(R+\rho)$ and $N_\eps$ are finite, the assertion follows from 
the inequality\\[-2mm]

 $\hspace{2cm} \#\{k\in\ZZ\;|\,|k|\le t/l(\gax)\}\le
2t/l(\gax)$. \qed\\[-2mm]
\begin{cor}\label{numinequiclass}
There exist  constants $b>0$ and $t_0>0$ 
\st for any $\gax\in\Gamma'$ primitive and all $\,t>t_0$ 
\begin{eqnarray*}
\#\{\gamma=\varphi \gax^k\varphi^{-1}\;|\,\varphi\in\Gamma,\;k\in\ZZ,\;
\piF(\varphi\gax^-)\in U_1,\,\piF(\varphi\gax^+)\in U_2,\,\ && \\
 l(\gamma)\le t \} & \le &  b \cdot t^r\,.
\end{eqnarray*}
\end{cor}
\prf\ We use the notation from the previous  lemma and notice that by choice
of $c_0>0$ the conditions 
$\piF(\varphi\gax^-)\in U_1$ and $\piF(\varphi\gax^+)\in U_2$ imply that   
\mbox{$\varphi\Ax(\gax)\cap B_\xo(c_0)\ne \emptyset$.}

Since $\overline{B_\xo(c_0)}\subset\XX$ is compact, there exist finitely many points $\{x_1, x_2,\ldots, x_m\}\subseteq \overline{B_\xo(c_0)}$ \st the balls $B_{x_i}(\rho)$, $1\le i\le m$, cover  $ \overline{B_\xo(c_0)}$.
Hence if $\piF(\varphi\gax^-)\in U_1$ and $\piF(\varphi\gax^+)\in U_2$, there exists
$j\in \{1,2,\ldots, m\}$ \st 
$\varphi\Ax(\gax)\cap\big(B_{x_j}(\rho)\cap
B_\xo(c_0)\big)\ne \emptyset$. We conclude\\[-1mm]

$\hspace{-5mm} \#\{\gamma=\varphi
\gax^k\varphi^{-1}\;|\,\varphi\in\Gamma,\;k\in\ZZ,\;
\piF(\varphi\gax^-)\in U_1,\,\piF(\varphi\gax^+)\in U_2,\, l(\gamma)\le t
\}\le m\cdot a\cdot t^r\,.$\qed\\


\section{The proof of the main theorem}

We will now state two more lemmata in order to obtain the upper bound and to relate $P(t)\,$ to
$\,N_\Gamma(R;A,B)$ for appropriate sets $A,B\subset K/M$. 
We use the notation from sections~\ref{Schottky} and \ref{lowbound}
and fix 
$$\eps <\min\{\; \delta,\; \eps_\Gamma(\{\piF(\gax_1^-)\},U_1),\;  \eps_\Gamma(\{\piF(\gax_1^+)\},U_2)\;\}\,.$$
\begin{lem}\label{fixedpointposition}
There exists $T>0$ 
\st for all $\gamma\in\Gamma$ with 
$$ d(\xo,\gamo)\ge T\,,\ \angle_\xo(\gamo, \piF(\gax_1^+))<\eps/6\,\ \
\mbox{and}\quad  \angle_\xo(\gamma^{-1}\xo, \piF(\gax_1^-))<\eps/6$$ 
either $g:=\gamma$ or $g:=\gamma^{-1}$ satisfies $\piF(g^-)\in U_1$ and $\piF(g^+)\in U_2$.
\end{lem} 
\prf\ 
If $\gamma\in\Gamma$ satisfies $\angle_\xo(\gamo,\piF(\gax_1^+))<\eps/6$
and $\angle_\xo(\gamma^{-1}\xo,\piF(\gax_1^-))<\eps/6$, then there exist
$\xi^-$, $\xi^+\in\regrand$ with $\piF(\xi^-)=\piF(\gax_1^-)$,
$\piF(\xi^+)=\piF(\gax_1^+)$, $\angle_\xo(\gamma^{-1}\xo,\xi^-)<\eps/4$ and
$\angle_\xo(\gamma\xo,\xi^+)<\eps/4$. Moreover, there exists
$R>0$ with the property that for every such $\gamma$ with $d(\xo,\gamo)\ge R $, and
for all $x\in \overline{B_\xo(c_0)}$ we have
$$\angle_x(\gamma\xo, \xi^+)<\eps/3\,,\ \angle_\xo(\gamma^{-1} x,
\xi^-)<\eps/3\,,\ \angle_x(\gamma^{-1}\xo, \xi^-)<\eps/3\,\ \an\quad \angle_\xo(\gamma x, \xi^+)<\eps/3\,.$$ 
Let $kM\in \overline{U_1}\subset U_1'$ arbitrary, and $\zeta\in\horinf(\xi^+)$ \st
$\piF(\zeta)=kM$. Denote $x\in\XX$ the orthogonal projection of $\xo$ to
the unique flat containing both $\zeta$ and $\xi^+$ in its boundary. 
Then $d(\xo,x)<c_0$ by definition of $c_0$, and we have
$\angle_x(\gamma \xo,\zeta)=\pi-\angle_x(\gamma\xo,\xi^+)>\pi-\eps/3$. From
$\angle_{\gamma \xo}(\zeta,x)+\angle_{x}(\gamma\xo,\zeta)\le \pi\ $
we therefore obtain
$$\angle_\xo(\gamma^{-1}\zeta,\gamma^{-1} x)=\angle_{\gamma\xo}(\zeta,x)\le \pi-\angle_x(\gamma\xo,\zeta)<\eps/3 \,,$$ 
hence
$\angle_\xo(\gamma^{-1}\zeta,\xi^-)\le \angle_\xo(\gamma^{-1}\zeta,\gamma^{-1}
x)+\angle_\xo(\gamma^{-1} x,\xi^-)<2\eps/3$. 

We now let $T\ge R$ \st every $\gamma\in\Gamma$ with $d(\xo,\gamma\xo)\ge T$ satisfies $\angle(H(\xo,\gamma\xo),P_\Gamma)<\eps/12$. We write $\xi^-=(k^-,H)$, where $k^-M=\piF(\gax_1^-)$ and $H\in\aL^+_1$ is the direction of $\xi^-$. If $H\in P_\Gamma$, then $\gamma^{-1}\zeta\in H_\Gamma^{\eps}(\{\piF(\gax_1^-)\})$, hence by choice of $\eps$ we have $\piF(\gamma^{-1}\zeta)=\gamma^{-1} kM\in U_1$. 

If $H\notin P_\Gamma$ we choose $H'\in P_\Gamma$ \st $\angle (H(\xo,\gamma^{-1}\xo),P_\Gamma)=\angle (H(\xo,\gamma^{-1}\xo),H')$. Then for $\eta^-:=(k^-,H')$ and $\gamma\in\Gamma$ with $d(\xo,\gamma\xo)\ge T$ we estimate 
\begin{eqnarray*}
\angle_\xo(\xi^-,\eta^-) &=&\angle (H,H')\le \angle (H,H(\xo,\gamma^{-1}\xo))+\angle(H(\xo,\gamma^{-1}\xo),H')\\
&\le &  \angle_\xo(\xi^-,\gamma^{-1}\xo)+\angle (H(\xo,\gamma^{-1}\xo),P_\Gamma) <\frac{\eps}4 +\frac{\eps}{12}=\frac{\eps}3\,.
\end{eqnarray*}
We conclude that $\angle_\xo(\gamma^{-1}\zeta,\eta^-)\le \angle_\xo(\gamma^{-1}\zeta,\xi^-)+\angle_\xo(\xi^-,\eta^-)<2\eps/3 +\eps/3=\eps$, hence also $\piF(\gamma^{-1}\zeta)=\gamma^{-1} kM\in U_1$. This proves $\gamma^{-1} \overline{U_1}\subseteq U_1$. 

Analogously we can show that $\gamma \overline{U_2}\subseteq U_2$, hence  
from Brouwer's
fixed point theorem it follows that  $\gamma$ possesses fixed
points in $U_1$ and $U_2$. Since every element in $\Gamma$ is regular axial by
Proposition~\ref{LimSetDescription}, the attractive and repulsive fixed point
$\gamma^+$, $\gamma^-$ are contained in the regular boundary. The assertion
now follows from 
the fact that\\[-2mm]

 $\hspace{3cm} \Fix^B(\gamma)\cap \bigcup_{i=1}^{2l} U_i =\{\piF(\gamma^+),\piF(\gamma^-)\}\,.$\qed\\
 

The following lemma will be crucial for the upper bound in our main theorem. It holds for any finitely generated free and discrete isometry group of a Hadamard manifold, hence in particular for our Schottky groups acting on a symmetric space of arbitrary rank. Notice that if $\Gamma=\langle \gax_1,\gax_2,\ldots, \gax_l\rangle$ is a free group, then every element $\gamma\in\Gamma$ can be written uniquely as  
$\gamma=s_1s_2\ldots s_n$ with letters $s_i\in S:=\{\gax_1,\gax_1^{-1},\ldots, \gax_l,\gax_l^{-1}\}$ \st $s_{i+1}\ne (s_i)^{-1}$ for $1\le i\le n-1$. We say that $\gamma$ is {\hl very reduced}, if $s_{n}\ne (s_1)^{-1}$. 
\begin{lem}\label{numberofconjugates}
Put $d:=\max \{ d(\xo,\gax_i\xo)\;|\, 1\le i\le l\}$. Then for any
$\gamma\in\Gamma$ 
there exist at least $[l(\gamma)/d-2]$ different very reduced conjugates of
$\gamma$.  
\end{lem}
\prf\ We write $\gamma\in\Gamma$ as above, i.e. $\gamma=s_1's_2'\cdots s_m'$
with $s_i'\in S$  and $s_{i+1}'\ne (s_i')^{-1}$ for $1\le i\le m-1$. If $s_1'=(s_{m}')^{-1}$, we conjugate $\gamma$ by $(s_1')^{-1}$ and get $\gamma'=s_2'
s_3'\cdots s_{m-1}'$. Repeating this procedure as long as possible, we obtain a
conjugate $\gamma(1)=s_1s_2\cdots s_n$ of $\gamma$ \st $s_1\ne s_n^{-1}$,
$n\le m$.  

For $2\le j\le n$ we put $\gamma(j):=s_js_{j+1}\cdots s_ns_1\cdots s_{j-1}$.
Then $\gamma(j)$ is conjugate to $\gamma(1)$ by $u(j):=(s_1s_2\cdots
s_{j-1})^{-1}$, and $\gamma(j)$  is different from $\gamma(i)$ for $i\ne j$ because there are no relations in
$\Gamma$. Since all elements $\gamma(j)$, $1\le j\le n$, are very reduced by construction, there exist at least $n-1$ different very reduced conjugates of
$\gamma$. 

In order to relate $n$ to $l(\gamma)$, we remark that 
\begin{eqnarray*}
l(\gamma) &= & l(\gamma(1))\le d(\xo,\gamma(1)\xo)\le d(\xo,
s_1\xo)+d(s_1\xo,s_1s_2\xo)+\cdots+d(s_1s_2\cdots s_{n-1}\xo, \gamma(1)\xo)\\
&= & d(\xo,s_1\xo)+d(\xo, s_2\xo)+\cdots +d(\xo,s_n\xo)\le  n\cdot \max_{s\in
  S} d(\xo,s\xo)=n\cdot d\,,\end{eqnarray*}
hence $n\ge l(\gamma)/d$.\qed\\

The following statement due to Y.~Benoist will finally give the upper bound. Recall Definition~\ref{thetavec} and the definition of the translation vector of an axial isometry from the beginning of section~\ref{ConjClassClosed}.
\begin{prp}(\cite[Section~4.1]{MR1437472})\\
Let $\Gamma<\is^o(\XX)$ be a Schottky group of a globally symmetric space $\XX$ as described in section~\ref{Schottky}. Then there exists a constant $M\ge 0$ \st every very reduced element $\gamma\in\Gamma$ satisfies $\Vert L(\gamma)-H(\xo,\gamma\xo)\Vert\le M$. 
\end{prp} 
Since $\Vert H(\xo,\gamma\xo)\Vert =d(\xo,\gamma\xo)$ and $\Vert L(\gamma)\Vert=l(\gamma)$, this proposition implies in particular that every very reduced element $\gamma\in\Gamma$ satisfies
\begin{equation}\label{LengthDistanceRelation} 
l(\gamma)\le d(\xo,\gamma\xo)\le \Vert H(\xo,\gamma\xo)-L(\gamma)+L(\gamma)\Vert\le M+ l(\gamma)\,.
\end{equation}
\begin{thr}\label{main}
Let $M$ be a locally symmetric space of noncompact type and rank $r\ge 1$  
with universal Riemannian covering manifold $\XX$ and fundamental group
isomorphic to a Schottky group $\Gamma\subset\is(\XX)$ as in section~\ref{Schottky}.
Then there exist constants $a>1$ and $t_0>0$ \st for all $t>t_0$
$$\frac{1}{a\cdot t^r}\cdot e^{\delta(\Gamma) t}\le P(t)\le \frac{a}{t}\cdot e^{\delta(\Gamma)t}\,.$$
\end{thr}
\prf\ For the proof of the upper bound, we let $t\ge 2+4d$ with $d>0$ as in Lemma~\ref{numberofconjugates} and consider a conjugacy class $[\gamma]$ with $t-1\le l([\gamma])< t$. By Lemma~\ref{numberofconjugates} there exist at least $\frac{t}{2d}$ different very reduced elements representing $[\gamma]$ with translation length equal to $l([\gamma])$, hence
$$ \Delta P(t):= P(t)-P(t-1)\le \frac{2d}{t}\cdot \#\{\gamma\in\Gamma\;|\, \gamma\ \mbox{very reduced}\,,\ l(\gamma)\in [t-1,t)\}\,.$$
Moreover, the estimate~(\ref{LengthDistanceRelation}) implies
$$ \#\{\gamma\in\Gamma\;|\, \gamma\ \mbox{very reduced}\,,\ l(\gamma)\in [t-1,t)\}\le N_\Gamma(t+M)\le b  e^{\delta(\Gamma)(t+M)}\,$$
by~(\ref{quintresult}). We summarize 
$$\Delta P(t)\le \frac{b'}{t} e^{\delta(\Gamma)t}\qquad\mbox{for some constant}\  b'>0\,.$$ 
Now Lemma~3.2 in \cite{MR2083973}
states that for all $n\in\NN$ and $a>0$ 
$$\sum_{k=1}^n \frac{e^{ak}}{k}\le \mbox{const} \cdot \frac{e^{an}}{n}\,.$$ 
So for $n\in\NN$
with $n\ge 2+4d$ we have 
$$ P(n)=\sum_{k=1}^n \Delta P(k)\le b' \sum_{k=1}^n
\frac{e^{\delta(\Gamma)k}}{k}\le \mbox{const} \cdot \frac{e^{\delta(\Gamma)n}}{n}\,,$$
which implies the upper bound. 

If the rank of $\XX$ equals one, then a Schottky group $\Gamma\subset\is(\XX)$
is convex cocompact, hence there exists a compact set $W\subset\XX$ \st every
closed geodesic in $M$ possesses a lift which intersects $W$
nontrivially. Moreover, $\Gamma$ satisfies the necessary conditions for Theorem~1 in
\cite{Link2005Asymptotic}, so the lower bound for $P(t)$ follows directly from
this theorem using the fact that $P(t)=P(t;W)$. 

If $r=\rank(\XX)>1$, we let $U_1, U_2\subset K/M$ be the neighborhoods of $\piF(\gax_1^-)$,
$\piF(\gax_1^+)$ as in section~\ref{Schottky}. Then by Lemma~\ref{fixedpointposition}, there exist $\eps>0$ and $T>0$ \st 
\be && \#\{\gamma\in\Gamma\;|\, d(o,\gamo)\ge T,\ \angle_\xo(\gamma\xo, \piF(\gax_1^+))<\eps,\ \angle_\xo(\gamma^{-1}\xo, \piF(\gax_1^-))<\eps \}\\
&& \ \ \le\, \#\{\gamma\in\Gamma\;|\, \piF(\gamma^-)\in U_1,\;
\piF(\gamma^+)\in U_2\;\}=\#\Gamma' \,.\ee 
Hence using the notation of
(\ref{NeighborhoodDef}),  putting
$$ A:=\piF(H_\Gamma^\eps(\{\piF(\gax_1^+)\}))\subset U_2\,,\qquad B:=\piF(H_\Gamma^\eps(\{\piF(\gax_1^-)\}))\subset U_1\,,$$
and applying the
inequality $l(\gamma)\le d(\xo,\gamo)$, $\gamma\in\Gamma$, we obtain
$$\#\{\gamma\in\Gamma'\;|\, l(\gamma)\le t \}\ge N_\Gamma(t;A,B)-N_\Gamma(T;A,B)\,.$$
Using Corollary~\ref{numinequiclass} and Proposition~\ref{countwithsets} we conclude that for  $t> T_0:=\max\{T,R_0\}$  
\be P(t)&\ge & \#\{[\gamma]\;|\, \gamma\in\Gamma' \ \,
\mbox{with}\ \, l([\gamma])\le t\} \ge  \frac{1}{b\cdot t^r}\cdot
\#\{\gamma\in\Gamma'|\, l( \gamma)\le t\}\\
&\ge & \frac1{b\cdot t^r}\big(\frac{1}{a} \cdot e^{\delta(\Gamma) t}  - a\cdot
e^{\delta(\Gamma)T_0}\big)\,,
\ee 
which proves that the lower bound holds for $t$ sufficiently large.  \qed\\

\bibliography{References}

\vspace{1.5cm}
\noindent Gabriele Link\\
Mathematisches Institut II\\
Universit{\"a}t Karlsruhe\\
Englerstr.~2\\
76 128 Karlsruhe\\
e-mail:\ gabriele.link@math.uni-karlsruhe.de

\end{document}